\iffalse
\documentclass{mcom-l}
\usepackage{amssymb}
\usepackage{graphicx}
\usepackage{booktabs}

\newtheorem{cor}{Corollary}
\usepackage{url}
\newcommand{\li}{{\rm li}}
\usepackage[ruled,vlined,linesnumbered]{algorithm2e}
\usepackage{comment}
\newtheorem{thm}{Theorem}
\newtheorem{lem}{Lemma}

\begin{document}

\title[On the first sign change of $\theta(x) -x$]{On the first sign change of $\theta(x) -x$}

\author{D. J. Platt}
\address{Heilbronn Institute for Mathematical Research University of Bristol, Bristol, UK}
\email{dave.platt@bris.ac.uk}
\curraddr{}
\thanks{}

\author{T. S.  Trudgian}
\address{Mathematical Sciences Institute, The Australian National University, ACT 0200, Australia}
\email{timothy.trudgian@anu.edu.au}
\curraddr{}
\thanks{Supported by Australian Research Council DECRA Grant DE120100173}



\dedicatory{}

\begin{abstract}
\noindent 
Let $\theta(x) = \sum_{p\leq x} \log p$. We show that $\theta(x)<x$ for $2<x< 1.39\cdot 10^{17}$. We also show that there is an $x<\exp(727.951332668)$ for which $\theta(x) >x.$
\end{abstract}

\maketitle

\else
\documentclass[11pt]{article}
\usepackage{a4wide}

\usepackage{booktabs}
\usepackage[ruled,vlined,linesnumbered]{algorithm2e}

\usepackage{amsthm}
\usepackage{amsmath}
\usepackage{amssymb}
\newtheorem{thm}{Theorem}
\newtheorem{lem}{Lemma}
\usepackage{graphicx}
\usepackage{booktabs}
\newtheorem{Rem}{Remark}
\newtheorem{cor}{Corollary}
\usepackage{url}
\newcommand{\li}{{\rm li}}

\title{On the first sign change of $\theta(x) -x$}

\author{
D.J. Platt\\ Heilbronn Institute for Mathematical Research \\ University of Bristol, Bristol, UK\\ dave.platt@bris.ac.uk
\and
  Tim Trudgian\footnote{Supported by Australian Research Council DECRA Grant DE120100173.} \\
  Mathematical Sciences Institute \\
  The Australian National University, ACT 0200, Australia \\
  timothy.trudgian@anu.edu.au
}

\begin{document}

\maketitle

\begin{abstract}
\noindent 
Let $\theta(x) = \sum_{p\leq x} \log p$. We show that $\theta(x)<x$ for $2<x< 1.39\cdot 10^{17}$. We also show that there is an $x<\exp(727.951332668)$ for which $\theta(x) >x.$
\end{abstract}


\fi

\section{Introduction}
Let $\pi(x)$ denote the number of primes not exceeding $x$. The prime number theorem is the statement that 
\begin{equation}\label{put}
\pi(x) \sim \textrm{li}(x) = \int_{2}^{x} \frac{dt}{\log t}.
\end{equation}
One often deals not with $\pi(x)$ but with the less obstinate Chebyshev functions $\theta(x)=\sum_{p\leq x} \log p$ and $\psi(x) = \sum_{p^{m}\leq x} \log p$.
The relation (\ref{put}) is equivalent to
\begin{equation*}\label{put2}
\psi(x) \sim x, \quad \textrm{and } \quad \theta(x) \sim x.
\end{equation*} 

Littlewood \cite{LittlePi}, showed that $\pi(x) - \textrm{li}(x)$ and $\psi(x) - x$ change sign infinitely often. Indeed, (see, e.g.,  \cite[Thms 34 \& 35]{Inghambook}) he showed more than this, namely that
\begin{equation}\label{pear}
\begin{split}
\pi(x) - \textrm{li}(x) &= \Omega_{\pm} \left(\frac{x^{\frac{1}{2}}}{\log x}\log\log\log x\right),\\
\psi(x) - x &= \Omega_{\pm} (x^{\frac{1}{2}} \log\log\log x).
\end{split}
\end{equation}
By \cite[(3.36)]{RS1} we have
\begin{equation}\label{ing1}
\psi(x) - \theta(x) \leq 1.427\sqrt{x} \quad (x>1),
\end{equation}
which, together with the second relation in (\ref{pear}), shows that $\theta(x) - x$ changes sign infinitely often.

Littlewood's proof that $\pi(x) - \textrm{li}(x)$ changes sign infinitely often was ineffective: the proof did not furnish a number $x_{0}$ such that one could guarantee that $\pi(x) - \textrm{li}(x)$ changes sign for some $x\leq x_{0}$. Skewes \cite{SkewesII} made Littlewood's theorem effective; the best known result is that there must be a sign change less that $1.3971\cdot 10^{316}$ \cite{SDT}. On the other hand Kotnik \cite{Kotnik} showed that $\pi(x) < \textrm{li}(x)$ for all $2<x\leq 10^{14}$.

We turn now to the question of sign changes in $\psi(x) -x$ and $\theta(x) - x$.
There is nothing of much interest to be said about the first sign changes of $\psi(x)$:  for $x\in[0, 100]$ there are 24 sign changes. The problem of determining an interval in which $\psi(x) - x$ changes sign is much more interesting (as examined in \cite{MontgomeryVorhauer}) but it is not something we consider here.
As for sign changes in $\theta(x)$: Schoenfeld, \cite[p.\ 360]{Schoenfeld} showed that $\theta(x) < x$ for all $0<x\leq 10^{11}$. This range appears to have been improved by Dusart, \cite[p.\ 4]{Dusart} to $0<x\leq 8\cdot 10^{11}$. We increase this in

\begin{thm}\label{bound}
For $0<x\leq1.39\cdot 10^{17}$, $\theta(x) < x$.
\end{thm}
A result of Rosser \cite[Lemma 4]{Rossers} is
\begin{lem}[Rosser]
If $\theta(x) <x$ for $e^{2.4}\leq x \leq K$ for some $K$, then $\pi(x) < \textrm{li}(x)$ for $e^{2.4} \leq x \leq K$.
\end{lem}
This enables us to extend Kotnik's result by proving
\begin{cor}
$\pi(x) < \textrm{li}(x)$ for all $2<x \leq 1.39\cdot 10^{17}$.
\end{cor}
 Rosser and Schoenfeld \cite[(3.38)]{RS1}, proved
\begin{equation}\label{rs1}
\psi(x) - \theta(x) - \theta(x^{\frac{1}{2}}) < 3x^{\frac{1}{3}}, \quad (x>0).
\end{equation}
Table 3 in \cite{Faber} gives us the bound $|\psi(x) -x| \leq 7.5\cdot 10^{-7}x$, which is valid for all $x\geq e^{35} > 1.5 \cdot 10^{15}$. This, together with  (\ref{rs1}) and Theorem \ref{bound}, enables us to make the following improvement to two results of Schoenfeld \cite[(5.1*) and (5.3*)]{Schoenfeld}. 
\begin{cor}
For $x>0$
\begin{equation*}\label{golf}
\theta(x) < (1 + 7.5\cdot 10^{-7})x, \quad \psi(x) - \theta(x) < (1 + 7.5\cdot 10^{-7})\sqrt{x} + 3x^{\frac{1}{3}}.
\end{equation*}
\end{cor}

We now turn to the question of sign changes in $\theta(x) -x$. In \S \ref{subsec:crossover} we prove
\begin{thm}\label{sign}
There is some $x\in[\exp(727.951332642),\exp(727.951332668)]$ for which $\theta(x)>x$.
\end{thm}
 
Throughout this article we make use of the following notation. For functions $f(x)$ and $g(x)$ we say that $f(x) =\mathcal{O}^{*}(g(x))$ if $|f(x)| \leq g(x)$ for the range of $x$ under consideration.

\section{Outline of argument}
The explicit formula for $\psi(x)$ is \cite[p.\ 101]{Inghambook}
\begin{equation}\label{ef}
\psi_{0}(x) = \frac{\psi(x+0) + \psi(x-0)}{2} = x - \sum_{\rho} \frac{x^{\rho}}{\rho} - \frac{\zeta'}{\zeta}(0) - \frac{1}{2} \log\left( 1 - \frac{1}{x^{2}}\right).
\end{equation}
Since 
$$\psi(x) = \theta(x) + \theta(x^{\frac{1}{2}}) + \theta(x^{\frac{1}{3}}) + \ldots,$$
we can manufacture an explicit formula for $\theta(x)$. Using (\ref{rs1}) and (\ref{ef}) we find that
\begin{equation}\label{tex}
\theta(x) -x > -\theta\left(x^{\frac{1}{2}}\right) - \sum_{\rho} \frac{x^{\rho}}{\rho} - \frac{\zeta'}{\zeta}(0) - 3 x^{\frac{1}{3}}.
\end{equation}
One can see why $\theta(x)<x$ `should' happen often. On the Riemann hypothesis $\rho = \frac{1}{2} + i\gamma$; since $\gamma \geq 14$ one expects the dominant term on the right-side of (\ref{tex}) to be $-\theta\left(x^{\frac{1}{2}}\right)$.

We proceed in a manner similar to that in Lehman \cite{LehmanPiLi}.
Let $\alpha$ be a positive number. We shall make frequent use of the Gaussian kernel $K(y) = \sqrt{\frac{\alpha}{2\pi}} \exp(-\frac{1}{2} \alpha y^{2})$, which has the property that $\int_{-\infty}^{\infty} K(y)\, dy = 1$.

Divide both sides of (\ref{tex}) by $x^{\frac{1}{2}}$, make the substitution $x\mapsto e^{u}$ and integrate against $K(u-\omega)$. This gives
\begin{equation}\label{main}
\begin{split}
\int_{\omega-\eta}^{\omega + \eta}&K(u-\omega)e^{\frac{u}{2}}\left\{ \theta(e^{u}) - e^{u}\right\}\, du > -\int_{\omega - \eta}^{\omega + \eta} K(u-\omega)\theta\left(e^{\frac{u}{2}}\right)e^{-\frac{u}{2}}\, du\\
& - \sum_{\rho} \frac{1}{\rho} \int_{\omega - \eta}^{\omega + \eta} K(u-\omega) e^{u(\rho - \frac{1}{2})}\, du - \frac{\zeta'(0)}{\zeta(0)}\int_{\omega - \eta}^{\omega + \eta}K(u-\omega)e^{-\frac{u}{2}}\, du \\
&- 3\int_{\omega - \eta}^{\omega + \eta}K(u-\omega) e^{-\frac{u}{6}}\, du  = -I_{1} - I_{2} - I_{3} - I_{4},
\end{split}
\end{equation}
say. The interchange of summation and integration may be justified by noting that the sum over the zeroes of $\zeta(s)$ in (\ref{tex}) converges boundedly in $u\in[\omega - \eta, \omega + \eta]$.
Noting that $\zeta'(0)/\zeta(0) = \log 2\pi$, we proceed to estimate $I_{3}$ and $I_{4}$ trivially to show that
\begin{equation*}
0 < I_3 < e^{-\frac{\omega-\eta}{2}} \log 2\pi, \quad
0 < I_4 < 3e^{-\frac{\omega-\eta}{6}}.
\end{equation*}
It will be shown in \S \ref{comp} that the contributions of $I_{3}$ and $I_{4}$ to (\ref{main}) are negligible --- this justifies our cavalier approach to their approximation.

We now turn to $I_{2}$. Let $A$ be the height to which the Riemann hypothesis has been verified, and let $T\leq A$ be the height to which we can reasonably compute zeroes to a high degree of accuracy --- we make this notion precise in \S \ref{comp}. Write $I_2=S_1+S_2$, where
\begin{equation*}
S_1=\sum\limits_{|\gamma| \leq A}\frac{1}{\rho}\int\limits_{\omega-\eta}^{\omega+\eta}K(u-\omega)e^{i\gamma u}\, du, \quad
S_2=\sum\limits_{|\gamma| > A}\frac{1}{\rho}\int\limits_{\omega-\eta}^{\omega+\eta}K(u-\omega)e^{(\rho-\frac{1}{2})u}\, du.
\end{equation*}
Our $S_{1}$ is the same as that used by Lehman in  \cite[pp.\ 402-403]{LehmanPiLi}. Using (4.8) and (4.9) of \cite{LehmanPiLi} shows that
\begin{equation*}
S_1=\sum\limits_{|\gamma|\leq T} \frac{e^{i\gamma\omega}}{\rho}e^{-\gamma^2/2\alpha} +E_{1},
\end{equation*}
where
\begin{equation*}\label{berry}
|E_1|<0.08\sqrt{\alpha}e^{-\alpha\eta^2/2}+e^{-T^2/2\alpha}\left\{\frac{\alpha}{\pi T^2}\log\frac{T}{2\pi}+8\frac{\log T}{T}+\frac{4\alpha}{T^3}\right\}.
\end{equation*}

Lehman considers 
$$f_{\rho}(s) = \rho s e^{-\rho s} \li(e^{\rho s}) e^{-\alpha(s-w)^{2} /2},$$
whence we writes his analogous version of $S_{2}$ as a function of $f_{\rho}(s)$ and then estimates this using integration by parts, Cauchy's theorem, and the bound
\begin{equation}\label{keys}
|f_{\rho}(s)| \leq 2 \exp(-\tfrac{1}{2} \alpha(s-w)^{2}).
\end{equation}
We consider the simpler function $f_{\rho}(s) = \exp(-\frac{1}{2}\alpha (s-w)^{2})$, which clearly satisfies (\ref{keys}). We may proceed as in \S 5 of \cite{LehmanPiLi} to deduce that
\begin{equation*}\label{S2}
|S_{2}| \leq A \log A e^{-A^{2}/(2a) + (w+ \eta)/2} \left\{ 4 \alpha^{-\frac{1}{2}} + 15 \eta\right\},
\end{equation*}
provided that 
\begin{equation*}\label{cons}
4A/w \leq \alpha \leq A^{2}, \quad 2A/\alpha \leq \eta < w/2.
\end{equation*}

All that remains is for us to estimate
\begin{equation*}
I_1=\int\limits_{\omega-\eta}^{\omega+\eta}\theta\left(e^{\frac{u}{2}}\right)e^{-\frac{u}{2}}K(u-\omega)\, du.
\end{equation*}
Table 3 in  \cite{Faber} and (\ref{ing1}) give us
\begin{equation}\label{thee}
|\theta(x) - x| \leq 1.5423\cdot 10^{-9} x, \quad x\geq e^{200},
\end{equation}
which gives
\begin{equation*}\label{ides}
I_1 < 1 + 1.5423\cdot 10^{-9}, \quad (\omega-\eta)\geq 400. 
\end{equation*}
Thus, we have
\begin{thm}\label{watch}
Let $A$ be the height to which the Riemann hypothesis has been verified, and let $T$ satisfy $0< T \leq A$. Let $\alpha, \eta$ and $\omega$ be positive numbers for which $\omega-\eta \geq 400$ and for which
\begin{equation*}
4A/\omega \leq \alpha \leq A^{2}, \quad 2A/\alpha \leq \eta \leq \omega/2.
\end{equation*}
Define $K(y) = \sqrt{\alpha/ (2\pi)} \exp(-\frac{1}{2} \alpha y^{2})$ and
\begin{equation}\label{pen}
I(\omega, \eta) = \int_{\omega - \eta}^{\omega + \eta} K(u - \omega) e^{-u/2} \left\{ \theta(e^{u}) - e^{u}\right\} \, du.
\end{equation}
Then
\begin{equation}\label{eyes}
I(\omega, \eta) \geq -1 - \sum_{|\gamma| \leq T} \frac{e^{i \gamma \omega}}{\rho} e^{-\gamma^{2}/(2\alpha)} - R_{1} - R_{2} - R_{3} - R_{4},
\end{equation}
where 
\begin{equation*}\label{chair}
\begin{split}
R_{1} & = 1.5423\cdot 10^{-9}\\
R_{2} & = 0.08\sqrt{\alpha}e^{-\alpha\eta^2/2}+e^{-T^2/2\alpha}\left\{\frac{\alpha}{\pi T^2}\log\frac{T}{2\pi}+8\frac{\log T}{T}+\frac{4\alpha}{T^3}\right\} \\
R_{3} & = e^{-(\omega - \eta)/2} \log 2\pi + 3e^{-(\omega - \eta)/6}  \\
R_{4} & = A (\log A) e^{-A^{2}/(2a) + (w+ \eta)/2} \left\{ 4 \alpha^{-\frac{1}{2}} + 15 \eta\right\}. \\
\end{split}
\end{equation*}
\end{thm}
We note that if one were to assume the Riemann Hypothesis for $\zeta$, then the $R_4$ term could be reduced. This would give us greater freedom in our choice of $\alpha$---see \S \ref{subsec:otherpars}.

Approximations different from (\ref{thee}) are available. For example, one could use Lemma 1 in \cite{Trudgianprime} to obtain $|\theta(x) - x| \leq 0.0045x/(\log x)^{2}$. One could also restrict the conditions in Theorem \ref{watch} to $\omega - \eta \geq 600$ using the slightly improved results from \cite{Faber} that are applicable thereto. Neither of these improves significantly the bounds in Theorem \ref{sign}.

We now need to search for values of $\omega$, $\eta$, $A$, $T$ and $\alpha$ for which the right-side of (\ref{eyes}) is positive.

\section{Computations}\label{comp}

\subsection{Locating a crossover}\label{subsec:crossover}
Consider the sum $\Sigma_{1} = \sum_{|\gamma| \leq T} \frac{e^{i \gamma \omega}}{\rho}$. We wish to find values of $T$ and $\omega$ for which this sum is small, that is, close to $-1$; for such values the sum that appears in (\ref{eyes}) should also small. 
Bays and Hudson \cite{Bays}, when considering the problem of the first sign change of $\pi(x) - \textrm{li}(x)$, identified some values of $\omega$ for which $\Sigma_{1}$ is small. We investigated their values: $\omega=405, 412, 437, 599, 686$ and $728$.

For $\omega$ in this range, we have $R_{1} = 1.5423\cdot 10^{-9}$ so we endeavour to choose the parameters $A, T, \alpha$ and $\eta$ to make the other error terms comparable.

\subsubsection{Choosing $A$}
We chose to rely on the rigorous verification of RH for $A=3.0610046\cdot 10^{10}$ by the second author \cite{Platt}. This computation also produced a database of the zeros below this height computed to an absolute accuracy of $\pm 2^{-102}$ \cite{Bober2012}.

\subsubsection{Choosing $T$}
As already observed, we have sufficient zeros to set $T=A\approx 3\cdot 10^{10}$ but, since summing over the roughly $10^{11}$ zeros below this height is too computationally expensive, we settled for $T=6,970,346,000$ (about $2\cdot 10^{10}$ zeros). Even then, computing the sum using multiple precision interval arithmetic (see \S \ref{subsec:sum}) takes about $40$ hours on an $8$ core platform. 

\subsubsection{Choosing the other parameters}\label{subsec:otherpars}

To get the finest granularity on our search (i.e.\ to be able to detect narrow regions where $\theta(x)>x$) we aim at setting $\eta$ as small as possible. This in turn means setting $\alpha$ (which controls the width of the Gaussian) as large as possible. However, to ensure that $R_4$ is manageable, we need $A^2/(2\alpha)>\omega/2$ or $\alpha<A^2/\omega$. 
 A little experimentation led us to
\begin{equation*}\label{choice}
\alpha=1,153,308,722,614,227,968, \quad \eta=\frac{933831}{2^{44}},
\end{equation*}
both of which are exactly representable in IEEE double precision.

\subsubsection{Summing over the zeros}\label{subsec:sum}
Since
\begin{equation*}
\frac{\exp(i\gamma\omega)}{\frac{1}{2}+i\gamma}+\frac{\exp(-i\gamma\omega)}{\frac{1}{2}-i\gamma}=\frac{\cos(\gamma\omega)+2\gamma\sin(\gamma\omega)}{\frac{1}{4}+\gamma^2},
\end{equation*}
the dominant term in $\Sigma_{1}$ is roughly $2\sin(\gamma\omega)/\gamma$. Though one might expect a relative accuracy of $2^{-53}$ when computing this in double precision, the effect of reducing $\gamma\omega$ mod $2\pi$ degrades this to something like $2^{-17}$ when $\gamma=10^9$ and $\omega=400$. We are therefore forced into using multiple precision, even though that entails a performance penalty perhaps as high as a factor of $100$. To avoid the need to consider rounding and truncation errors at all, we use the MPFI \cite{Revol2005} multiple precision interval arithmetic package for all floating point computations. Making the change from scalar to interval arithmetic probably costs us another factor of $4$ in terms of performance.

\subsubsection{Results}

\begin{figure}[tbp]
\centering
\fbox{\includegraphics[width=0.7\linewidth, angle=270]{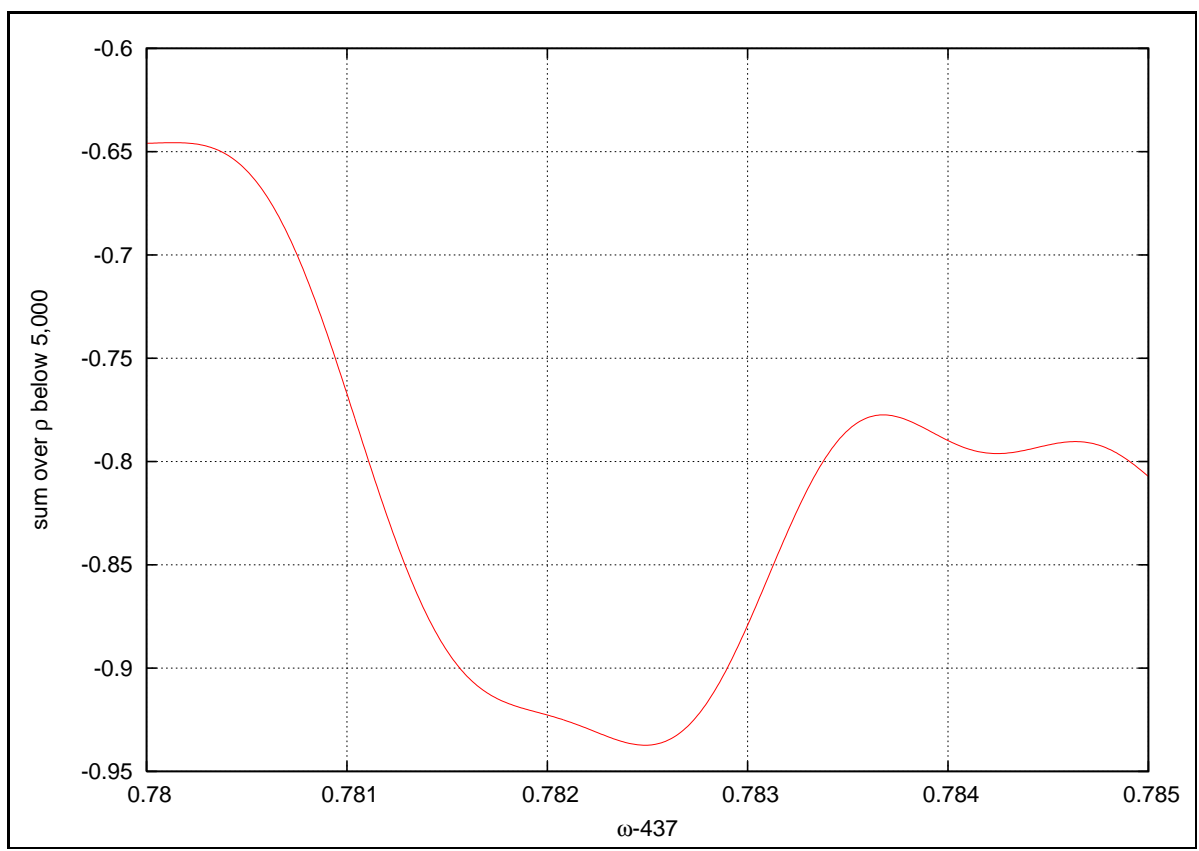}}
\caption{Plot of $\sum_{|\gamma| \leq 5000} \frac{e^{i\omega \gamma}}{\rho}$ for $\omega \in[437.78, 437.785]$.}
\label{fig:437}
\end{figure}

We initially searched the regions around $\omega=405, 412, 437, 599, 686$ and $728$ using only those zeros $\frac{1}{2} + i\gamma$ with  $0<\gamma<T=5,000$. Although these results were not rigorous, it was hoped that a sum approaching $-1$ would indicate a potential crossover worth investigating with full rigour. As an example, Figure \ref{fig:437} shows the results for a region near $\omega=437.7825$. This is some way from dipping below the $-1$ level and indeed a rigorous computation using the full set of zeros and with $\omega=437.78249$ fails to get over the line. The same pattern repeats for $\omega$ near $405, 412, 599 $ and $686$.

In contrast, we expected the region near $728$ to yield a point where $\theta(x)>x$. The lowest published interval containing an $x$ such that $\pi(x)>\li(x)$ is 
$$x\in[\exp(727.951335231),\exp(727.951335621)]$$ 
in \cite{SDT}. Since the error terms for $\theta(x)-x$ are tighter
than those for $\pi(x)-\li(x)$ this necessarily means that the same $x$ will satisfy $\theta(x)>x$. In fact, we can do better. Using $\omega=727.951332655$ we get
\begin{equation*}
\sum\limits_{|\gamma|\leq T} \frac{\exp(i\gamma\omega)}{\rho}\exp\left(-\frac{\gamma^2}{2\alpha}\right)\in[-1.0013360278,-1.0013360277].
\end{equation*}

We also have $R_{1}+ R_{2}+ R_{3} + R_{4}<1.7\cdot 10^{-9},$
so that
\begin{equation}\label{eq:b1}
\int_{\omega-\eta}^{\omega + \eta}K(u-\omega)e^{-u/2}\left\{ \theta(e^{u}) - e^{u}\right\}\, du > 0.0013360261.
\end{equation}

\subsubsection{Sharpening the Region}

Using the same argument as \cite[\S 9]{SDT}, we can analyse the tails of the integral (\ref{pen}) and sharpen the region considerably. Consider, for $\eta_0\in (0,\eta]$,
\begin{equation*}
T_{1}=\int\limits_{\omega+\eta_0}^{\omega+\eta}K(u-\omega)e^{-\frac{u}{2}}\left\{\theta\left(e^u\right)-e^u\right\} du,
\end{equation*}
and
\begin{equation*}
T_{2}=\int\limits_{\omega-\eta}^{\omega-\eta_0}K(u-\omega)e^{-\frac{u}{2}}\left\{\theta\left(e^u\right)-e^u\right\} du.
\end{equation*}
Another appeal to Table 3 in \cite{Faber}, and (\ref{ing1}), gives us 
\begin{equation*}
|\theta(x) - x| \leq 1.3082\cdot 10^{-9} x, \quad x\geq e^{700}.
\end{equation*}
Thus for $\omega-\eta>700$ we have
\begin{equation}\label{eq:tails}
\left|T_1\right|+\left|T_2\right|\leq 1.3082\cdot 10^{-9} (\eta-\eta_0)K(\eta_0)\left[e^{\frac{\omega+\eta}{2}}+e^{\frac{\omega-\eta_0}{2}}\right].
\end{equation}
Applying (\ref{eq:tails}) to (\ref{eq:b1}), we find we can take $\eta_0=\eta/4.2867$ so that
\begin{equation*}
\int_{\omega-\eta_0}^{\omega + \eta_0}K(u-\omega)e^{-u/2}\left\{ \theta(e^{u}) - e^{u}\right\}\, du > 2.75\cdot 10^{-6},
\end{equation*}
which proves Theorem \ref{sign}.
Therefore, there is at least one $u\in(\omega-\eta_{0}, \omega + \eta_{0})$ with $\theta(e^{u}) - e^{u}>0$. Owing to the positivity of the kernel $K(u-\omega)$ we deduce that there is at least one such $u$ with 
\begin{equation*}
\theta(e^{u}) - e^{u} > 2.75\cdot 10^{-6} e^{u/2} > 10^{152}.
\end{equation*}
Since $\theta(x)$ is non-decreasing this proves
\begin{cor}\label{cord}
There are more than $10^{152}$ successive integers $x$ satisfying $$x\in[\exp(727.951332642),\exp(727.951332668)],$$ for which $\theta(x)>x$.
\end{cor}
\subsection{A lower bound}

Having established an upper bound for the first time that $\theta(x)$ exceeds $x$, we now turn to a lower bound. A simple method would be to sieve all the primes $p$ less than some bound $B$, sum $\log p$ starting at $p=2$, and compare the running total each time to $p$. We set $B=1.39\cdot 10^{17}$ since this was required by the second author for another result in \cite{DudekPlatt}. By the prime number theorem we would expect to find about $3.5\cdot 10^{15}$ primes below this bound. Since this is far too many for a single thread computation we must look for some way of computing in parallel.

\subsubsection{A parallel algorithm}

We divide the range $[0,B]$ into contiguous segments. For each segment $S_j=[x_j,y_j]$ we set $T=\Delta=\Delta_{\textrm{min}}=0$. We look at the each prime $p_i$ in this segment, compute $l_i=\log p_i$, and add it to $T$. We set $\Delta=\Delta+l_i-p_i+p_{i-1}$ and $\Delta_{\textrm{min}}=\min(\Delta_{\textrm{min}},\Delta)$. Thus at any $p$, $\Delta_{\textrm{min}}$ is the maximum amount by which $\theta(p)$ has caught up with or gone further ahead of $p$ within this segment. After processing all the primes within a segment, we output $T$ and $\Delta_{\textrm{min}}$.

Now, for each segment $S_j=[x,y]$ the value of $\theta(x)$ is simply the sum of $T_k$ with $k<j$ and $\theta(y)=\theta(x)+T_j$. Furthermore, if $\theta(x)<x$ and $\theta(x)+\Delta_{\textrm{min}}>0$ then $\theta(w)<w$ for all $w\in[x,y]$.

\subsubsection{Results}

We implemented this algorithm in C++ using Kim Walisch's ``primesieve'' \cite{Walisch2012} to enumerate the primes efficiently, and the second author's double precision interval arithmetic package to manage rounding errors.

We split $B$ into $10,000$ segments of width $10^{13}$ followed by $390$ segments of width $10^{14}$. This pattern was chosen so that we could use Oliviera e Silva's tables of $\pi(x)$ \cite{Oliviera2012} as an independent check of the sieving process.

We used the $16$ core nodes of the University of Bristol Bluecrystal Phase III cluster \cite{ACRC2014} and we were able to utilise each core fully. In total we used about $78,000$ node hours. This established Theorem \ref{bound}.

\begin{figure}[tbp]
\centering
\fbox{\includegraphics[width=0.7\linewidth, angle=270]{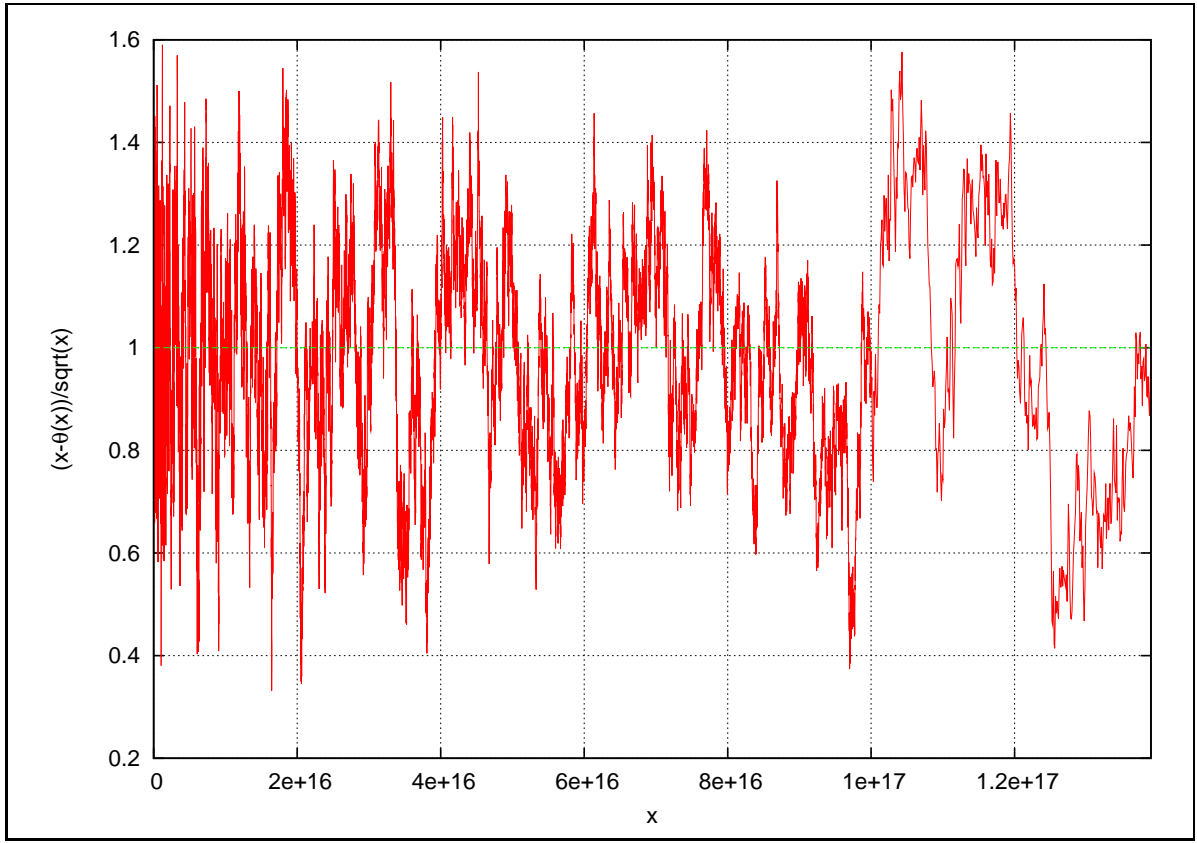}}
\caption{Plot of $\frac{x - \theta(x)}{\sqrt{x}}$.}
\label{fig:cheb}
\end{figure}

We plot $(x-\theta(x))/\sqrt{x}$ measured at the end of each segment in Figure \ref{fig:cheb}. As one would expect, this appears to be a random walk around the line $1$.

\bibliographystyle{plain}
\bibliography{themastercanada}

\end{document}